\documentclass[preprint,12pt]{amsart}

\usepackage[utf8]{inputenc}

\usepackage{graphicx}
\usepackage{graphicx}
\usepackage{amssymb}
\usepackage{enumerate}
\usepackage{amsthm}
\usepackage{textcomp}
\usepackage{graphicx,epsfig}
\usepackage{amsfonts,amssymb,enumerate}
\usepackage{latexsym,amsfonts,amssymb,amsmath,amscd,eufrak}
\usepackage{fancybox,shadow,color,graphicx,psfrag,float}
\usepackage{pstricks,pst-node,pst-text,pst-3d,pst-plot}
\usepackage{tikz,pgffor}

\newtheorem{theorem}{Theorem}[section]

\theoremstyle{definition}

\theoremstyle{remark}

\numberwithin{equation}{section}

\usepackage{amssymb}

\usepackage{amsmath}

\newcommand{\m}{\mathcal}

\title{A family of asymptotically bad wild towers of function fields}

 \author{M. Chara }  \address{ M. Chara: Researcher of CONICET at FIQ, Universidad Nacional del Litoral, Argentina, mchara@santafe-conicet.gov.ar}  
 \author{R. Toledano}
 \address { R. Toledano: Departamento de Matemática, FIQ. Universidad Nacional del Litoral, Argentina, ridatole@gmail.com} 
\thanks{This work was partially funded by Plan de Excelencia en Investigación Científica (PEIC I+D-2022-066)}
\begin{document}
\maketitle

\begin{abstract}

In \cite{badtowers} general conditions were given to prove the infiniteness of the genus of certain towers of function fields over a perfect field. It was shown that many examples where particular cases of those general results. In this paper the genus of a family of wild towers of function fields will be considered together with a result with less restrictive sufficient conditions for a wild tower to have infinite genus.
\end{abstract}

\section{Introduction} 
\label{intro}

The problem of the determination of whether or not a given sequence of function fields over a perfect field $K$ has finite  genus, has been of interest because the asymptotic behavior of sequences of codes constructed from asymptotically good towers of function fields has interesting applications in  coding theory and cryptography. In the books \cite{NX01}, \cite{GS07} and \cite{Stichbook09} can be found a good  variety of examples where this problem has been addressed for the case of recursive towers of function fields over a finite field. 

A  \emph{tower} of function fields over a perfect field $K$ is simply an infinite sequence $\m{T}=\{T_i\}_{i=0}^{\infty}$ of function fields $T_i$ over $K$ such that $T_{i}$ is a  subfield of $T_{i+1}$  for every $i\geq 0$. These mathematical objects turned out to be useful in Coding theory, Cryptography and related areas, as can be seen in \cite{NX01} and \cite{Stichbook09}.  For these applications more conditions are required on a tower, namely that  each extension $T_{i+1}/T_i$ be a finite and separable extension, $g(T_i)\rightarrow \infty$ as $i\rightarrow\infty$, where $g(T_i)$ denotes the genus of $T_i$, and also that $K$ be the full constant field of each $T_i$.

An important quantity  associated to a tower $\m{T}=\{T_i\}_{i=0}^{\infty}$ of function fields over $K$ is its \emph{genus} $\gamma(\mathcal{T})$ over $T_0$ which is defined as 
$$\gamma(\mathcal{T}):=\lim_{i \rightarrow \infty}\frac{g(T_i)}{[T_i:T_0]}.$$
It is proved in Chapter 7 of \cite{Stichbook09} that $0<\gamma(\mathcal{T})\leq \infty$ and for the above mentioned applications it is desirable to deal with towers of finite genus (see Chapter 7 of \cite{Stichbook09}). Thus it is interesting not only to have conditions ensuring the finiteness of the genus of a tower but also to be able to discard towers because of the infiniteness of their genus. 

We recall that {\it non skew recursive towers} of function fields are the most useful in the applications. Recursive means that the first function field in the tower is a rational function field, that is  $T_0=K(x_0)$ for some transcendental element $x_0$ over $K$, and all the others fields in the sequence can be defined recursively for every $i\geq 0$ as $T_{i+1}=T_{i}(x_{i+1})$, where $\{x_i\}_{i=0}^{\infty}$ is a sequence of transcendental elements over $K$ satisfying an equation of the form $$F(x_{i},x_{i+1})=0,$$ for some appropriate bivariate polynomial $F \in K[X,Y]$. A \emph{non-skew} recursive tower is a recursive tower in which the defining polynomial $F$ satisfies that $\deg_X(F)=\deg_Y(F)$.

Our main theoretical result in this paper is Theorem \ref{wildramif} where a general sufficient condition to have infinite genus in certain wild non-skew towers of function fields over a perfect field $K$ is given. A particular case of this condition was studied by Garcia and Stichtenoth in \cite{GS96}, and then generalized by Chara and Toledano in \cite{badtowers}. Other general conditions for the infiniteness of the genus are given in \cite{BeelenGS05}. We illustrate Theorem \ref{wildramif} with a new example of a family of non skew recursive wild towers having infinite genus.

Let $\m{T}=\{T_i\}_{i=0}^{\infty}$ be a tower of function fields over $K$. The tower $\mathcal{T}$ is said to be a \emph{tame} tower if each extension $T_{i+1}/T_i$ is tamely ramified, that is the ramification index $e(Q|P)$ is not divisible by $\textrm{Char}(K)$ for any place $Q$ of $T_{i+1}$ lying above a place $P$ of $T_i$. Otherwise the tower $\mathcal{T}$ is said to be a \emph{wild} tower, that is there is an extension $T_{i+1}/T_i$ and a  place $Q$ of $T_{i+1}$ lying above a place $P$ of $T_i$ such that the ramification index $e(Q|P)$ is divisible by $\textrm{Char}(K)$.

\section{Climbing a wild tower}\label{pyramid}

 In what follows we will use the symbol $d(Q|P)$ to denote the different exponent attached to a place $Q$ lying over a place $P$ (see, for example, Chapter 1 of \cite{NX01} for a quick review of the basic properties of the different exponent).

Also a  place defined by a monic and irreducible polynomial $f\in K[x]$ in a rational function field $K(x)$ will be denoted by $P_{f}$ or $P_{f(x)}$ if the occurrence of the transcendental element $x$ is needed. The symbol $P_\infty$ will be used to denote the place which is the only pole of $x$ in $K(x)$.

It was shown in \cite{badtowers} that in many cases the infiniteness of the genus of a wild tower $\m{T}$ can be proved by finding, for infinitely many fields $T_i$ in the tower,  a place $Q$ of $T_{i+1}$ lying over a place $P$ of $T_i$  such that $d(Q|P)\geq c_i\cdot [T_{i}:T_0]$ with $c_i\geq0$ and the series
\begin{equation}\label{seriesdiffexp}
\sum_{i=1}^{\infty}\frac{c_{i}}{[T_{i+1}:T_{i}]}
\end{equation} is divergent. We give in the next theorem  a general sufficient condition for the infiniteness of the genus of a recursive wild tower.

\begin{theorem}\label{wildramif}
Let $\m{T}=\{T_i\}_{i=0}^{\infty}$ be a non skew recursive tower of function fields over a perfect field $K$ of characteristic $p>0$ defined by a suitable bivariate polynomial $F\in K[x,y]$. Let $m=\deg_yF=\deg_xF$ and let us consider the basic function field $T=K(x,y)$ associated to  $\m{T}$. Suppose that there exist a monic and irreducible polynomial $f\in K[X]$ and a place $Q$ of $T$ lying above $P_{f(y)}$  (the zero of $f(y)$ in $K(y)$) and above $P_{f(x)}$ (the zero of $f(x)$ in $K(x)$) such that
\begin{enumerate}[(1)]
\item $m=e(Q|P_{f(y)})$ and $\gcd(m,p)=1$ (i.e. $P_{f(y)}$ is totally and tamely ramified in $K(x,y)$) and \label{1-examplewildramif}
\item  $\gcd(e(Q|P_{f(x)}),m)=1$.\label{2-examplewildramif}
\item  There exits a place $Q'$ of $K(x,y)$ lying above $P_{f(x)}$ such that $Q'|P_{f(x)}$ is wildly ramified and $\gcd(e(Q'|P_{f(x)}),m)=1$.
\end{enumerate}
Then $\gamma(\m{T})=\infty$.
 \end{theorem}
 \begin{proof}
Let $\{x_i\}_{i=0}^{\infty}$ be a sequence of transcendental elements over $K$ such that $T_0=K(x_0)$ and $T_{i+1}=T_i(x_{i+1})$ for $i\geq 0$. For $i\geq 1$ let us write $P_i=P_{f(x_{i})}$. By hypothesis we have a place $Q_i$ of $K(x_{i-1},x_i)$ lying above the places $P_{i-1}$ of $K(x_{i-1})$ and $P_i$ of $K(x_i)$ . Let $P'$ be a place of $T_{i+1}$  lying above $Q'$ (see Figure \ref{figu0} below). We show next that $P'$ lies above $P_j=P_{f(x_{j})}$ for $j=0,\ldots,i$. Since $Q'$ lies above $P_i$ we have that $P'$ lies above $P_i$. Thus $P'\cap K(x_{i-1},x_i)$ is a place of $K(x_{i-1},x_i)$ lying above $P_i$. By hypothesis $Q_i$ is the only place of $K(x_{i-1},x_i)$ lying above $P_i$. Then $P'\cap K(x_{i-1},x_i)=Q_i$ so that $P'$ lies above $P_{i-1}$. We have now that $P'\cap K(x_{i-2},x_{i-1})$ is a place of $K(x_{i-2},x_{i-1})$ lying above $P_{i-1}$. By hypothesis $Q_{i-1}$ is the only place of $K(x_{i-2},x_{i-1})$ lying above $P_{i-1}$. Then $P'\cap K(x_{i-2},x_{i-1})=Q_{i-1}$ so that $P'$ lies above $P_{i-2}$. Continuing in this way we obtain the desired property.

By  Abhyankar's lemma (see \cite[Theorem 3.9.1]{Stichbook09}) we have the following ramification picture where $n=e(Q|P_{f(x)})$, $r=e(Q'|P_{f(x)})$ and $P$ is the place of $T_i$ lying under $P'$. By \eqref{1-examplewildramif} we have that $P$ lies above $P_i$.

 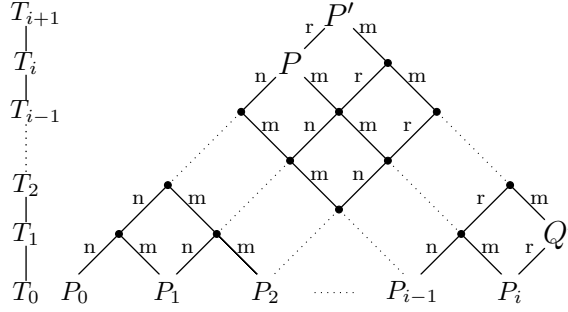
\begin{figure}[h!t]
        \begin{center}
  \begin{tikzpicture}[scale=0.65]
 \draw[line width=0.5 pt](-0.9,0)--(-0.9,2.2);
  \draw[line width=0.5 pt](-0.9,3.7)--(-0.9,5.7);
   \draw[line width=0.5 pt, dotted](-0.9,2.2)--(-0.9,3.7);
   \draw[white,  fill=white](-0.9,0) circle (0.23 cm);
      \draw[white,  fill=white](-0.9,1.2) circle (0.23 cm);
         \draw[white,  fill=white](-0.9,2.2) circle (0.23 cm);
            \draw[white,  fill=white](-0.9,3.7) circle (0.23 cm);
               \draw[white,  fill=white](-0.9,4.7) circle (0.23 cm);
                  \draw[white,  fill=white](-0.9,5.7) circle (0.23 cm);
\node at(-0.9,0){\footnotesize{$T_0$}};
\node at(-0.9,1.2){\footnotesize{$T_1$}};
\node at(-0.9,2.2){\footnotesize{$T_2$}};
\node at(-0.7,3.7){\footnotesize{$T_{i-1}$}};
\node at(-0.9,4.7){\footnotesize{$T_{i}$}};
\node at(-0.7,5.7){\footnotesize{$T_{i+1}$}};
      \draw[dotted](2,2.2)--(3.5,3.7);
     \draw[dotted](3,1.2)--(4.5,2.7);
     \draw[dotted](4,0.2)--(5.5,1.7);
     \draw[dotted](8,1.2)--(6.5,2.7);
     \draw[dotted](7,0.2)--(5.5,1.7);
\draw[dotted](5,0.0)--(5.9,0.0);
\draw[dotted](9,2.2)--(7.5,3.7);
    \draw[line width=0.5 pt](0,0.2)--(2,2.2)--(4,0.2)--(3,1.2)--(2,0.2)--(1,1.2);
  \draw[line width=0.5 pt](7,0.2)--(9,2.2);
\draw[line width=0.5 pt](10,1.2)--(9,2.2);
\draw[line width=0.5 pt](10,1.2)--(9,0.2)--(8,1.2);
  \draw[line width=0.5 pt](5.5,3.7)--(4.5,2.7)--(5.5,1.7)--(6.5,2.7)--(5.5,3.7)--(4.5,4.7)--(3.5,3.7)--(4.5,2.7);
 \draw[line width=0.5 pt](6.5,2.7)--(7.5,3.7)--(5.5,5.7)--(4.5,4.7);
  \draw[line width=0.5 pt](5.5,3.7)--(6.5,4.7);
   \draw[white,  fill=white](0,0) circle (0.4 cm);
    \draw[white,  fill=white](2,0) circle (0.4 cm);
  \draw[white,  fill=white](4,0) circle (0.4 cm);
 \draw[white,  fill=white](7,0) circle (0.4 cm);
  \draw[white,  fill=white](9,0) circle (0.4 cm);
   \draw[white,  fill=white](11,0) circle (0.4 cm);
  \node at(0.1,0){\footnotesize{$P_0$}};
   \node at(2,0){\footnotesize{$P_1$}};
    \node at(4,0){\footnotesize{$P_2$}};
      \node at(7,0){\footnotesize{$P_{i-1}$}};
     \node at(9,0){\footnotesize{$P_i$}};
  \draw[fill=black](1,1.2) circle (0.07 cm);
     \node at(1,1.2){$\,$};
       \draw[fill=black](3,1.2) circle (0.07 cm);
     \node at(3,1.2){$\,$};
       \draw[fill=black](2,2.2) circle (0.07 cm);
     \node at(2,2.2){$\,$};
     \draw[fill=black](5.5,1.7) circle (0.07 cm);
     \draw[fill=black](4.5,2.7) circle (0.07 cm);
     \draw[fill=black](3.5,3.7) circle (0.07 cm);
     \draw[fill=black](6.5,2.7) circle (0.07 cm);
     \draw[fill=black](5.5,3.7) circle (0.07 cm);
     \draw[white, fill=white](4.5,4.7) circle (0.4 cm);
    \draw[white, fill=white](5.5,5.7) circle (0.4 cm);
     \draw[fill=black](8,1.2) circle (0.07 cm);
     \draw[white, fill=white](10,1.2) circle (0.4 cm);
       \node at(4.5,4.7){$P$};
       \node at(5.5,5.7){$P'$};
         \node at(10,1.2){$Q'$};
  \draw[fill=black](6.5,4.7) circle (0.07 cm);
  \draw[fill=black](7.5,3.7) circle (0.07 cm);
  \draw[fill=black](9,2.2) circle (0.07 cm);
       \node at(0.4,0.9){\tiny{n}};
          \node at(1.6,0.9){\tiny{m}};
           \node at(2.4,0.9){\tiny{n}};
            \node at(3.6,0.9){\tiny{m}};
             \node at(7.4,0.9){\tiny{n}};
             \node at(8.6,0.9){\tiny{m}};
              \node at(9.4,0.9){\tiny{r}};
                   \node at(1.4,1.9){\tiny{n}};
          \node at(2.6,1.9){\tiny{m}};
              \node at(8.4,1.9){\tiny{r}};
               \node at(9.6,1.9){\tiny{m}};
                  \node at(5.1,2.4){\tiny{m}};
                    \node at(5.9,2.4){\tiny{n}};
                      \node at(6.9,3.4){\tiny{r}};
                       \node at(6.1,3.4){\tiny{m}};
                        \node at(4.9,3.4){\tiny{n}};
                        \node at(4.1,3.4){\tiny{m}};
                        \node at(3.9,4.4){\tiny{n}};
                         \node at(5.1,4.4){\tiny{m}};
                         \node at(5.9,4.4){\tiny{r}};
                           \node at(7.1,4.4){\tiny{m}};
                         \node at(6.1,5.4){\tiny{m}};
                           \node at(4.9,5.4){\tiny{r}};
 \end{tikzpicture}
  \caption{Ramification of $P_i$ in $F_i$ in Theorem \ref{wildramif}}\label{figu0}
\end{center}\end{figure}
Since $T_{i+1}=T_i\cdot K(x_i,x_{i+1})$ and $\gcd(m,p)=1$, the transitivity formula for the different exponent (see \cite[Corollary 3.4.12]{Stichbook09})
implies that
\begin{align*}
d(P'|P) &= d(P'|P_i)-e(P'|P)d(P|P_i)\\
        &= e(P'|Q')d(Q'|P_i)+d(P'|Q')-e(P'|P)d(P|P_i)\\
        &= m^id(Q'|P_i)+m^i-1-r(m^i-1)\\
        &\geq m^ir+(m^i-1)(1-r)> \frac{m^i}{2}+\frac{m^i}{2}(1-r)\\
        &= \frac{m^i}{2}=\frac{[T_i:T_0]}{2}\,,
\end{align*}
where in the last equality we have used that the tower $\m{T}$ non skew. We have that  the series \eqref{seriesdiffexp} is divergent with $c_i=2$ because $[T_{i+1}:T_i]=m$ for all $i\geq 0$.
 Therefore we conclude that $\gamma(\m{T})=\infty$ as desired.
 \end{proof}

\section{A new family of skew recursive wild towers with infinite genus}
Let $q$ be a prime $p$ power, $a\in K$ and let $g$ be a polynomial over $K$ such that $g(a)\neq 0$ with $\deg(g)< m=q+1$. Let us consider the sequence $\mathcal{T}=\{T_i\}_{i=0}^\infty$ of function fields over $K$ recursively defined by the equation
\[(y-a)^m+b(y-a)=\frac{(x-a)^m}{g(x)},\]
where $0\neq b\in K$ and suppose also that $\gcd(m-\deg(g),m)=1$

For simplicity we can assume that $K$ is algebraically closed. There is no harm with this assumption because we are going to prove that the genus of $\mathcal{T}$ goes to infinity, and the genus of any function field does not change for constant field extensions. 
 
First of all we have that the extension $K(x,y)/K(x)$ is of degree $m$. In fact, using Eisenstein's irreducibility criterion (see, for example, Proposition 3.1.15 of \cite{Stichbook09}) with the place $P_\infty$ we see that the polynomial $\varphi(T)=T^{q+1}-aT^q+(b-a^q)T+(a^q-b)a-(x-a)^m/g(x)$ over $K(x)$
 is irreducible over $K(x)$ and that $\varphi(T)$ is the minimal polynomial of $y$ over $K(x)$. Thus $K(x,y)/K(x)$ is an extension of degree $m$ and, since we are assuming that $K$ is algebraically closed, we also have that $K$ is the full constant field of $K(x,y)$.

  Let $Q$ be a zero of $y-a$ in $K(x,y)$. Since $\nu_Q(y-a)>0$, then $\nu_Q((y-a)^m)=m\nu_Q(y-a)>\nu_Q(y-a)$ since $m\geq 2$. Thus, from the Strict triangular inequality (Lemma 1.1.11 of \cite{Stichbook09}) we have that
\[\nu_Q((y-a)^{m}+b(y-a))=\nu_Q(b(y-a))=\nu_Q(y-a)>0.\]
Let $P$ be a place of $K(x)$ lying below $Q$.  Then
\[0<\nu_Q((y-a)^{m}+b(y-a))=e(Q|P)(m\nu_P(x-a)-\nu_P(g(x))),\]
so that $m\nu_P(x-a)>\nu_P(g(x))$. If $P=P_\infty$ then $\nu_P((x-a)^m/g(x))<0$, 
a contradiction. Thus $P\neq P_\infty$ and then $\nu_P(g(x))\geq 0$ because $g$ is a polynomial. Therefore $\nu_P(x-a)>0$ and we must have that $P=P_{x-a}$, the only zero of $x-a$ in $K(x)$. Thus $Q$ lies above $P_{x-a}$ and $P_{y-a}$, the only zero of $y-a$ in $K(y)$.
 
Now let $Q'$ be a zero of $(y-a)^q+b$ in $K(x,y)$. Then $\nu_{Q'}(y-a)=0$   so that
\[\nu_{Q'}((y-a)^{m}+b(y-a))=\nu_{Q'}((y-a)^q+b)>0.\]
Thus $Q'$ is also a zero of $(y-a)^{m}+b(y-a)$ in $K(x,y)$. Let $R$ be a place of $K(x)$ lying below $Q'$. Then
\[0<\nu_{Q'}((y-a)^{m}+b(y-a))=e(Q'|R)(m\nu_R(x-a)-\nu_R(g(x))),\]
so that $m\nu_R(x-a)>\nu_R(g(x))$. As before, must have that $R=P_{x-a}$. Notice that $\nu_{P_{x-a}}(g(x))=0$ because $g(x)$ and $x-a$ are coprime polynomials by hypothesis. Since $K$ is algebraically closed, there is $c\in K$ such that $c^q=b$. Then
$(y-a)^q+b=(y-a+c)^q$,
so that $q\nu_{Q'}(y-a+c)= \nu_{Q'}((y-a)^q+b)$. On the other hand, since $\nu_{Q'}(y-a)=0$ and $\nu_{P_{x-a}}(g(x))=0$, we have that
\begin{align*}
	\nu_{Q'}((y-a)^q+b)&=\nu_{Q'}((y-a)^{m}+b(y-a))\\
	&=e(Q'|P_{x-a})\nu_{P_{x-a}}((x-a)^m/g(x))=m\,e(Q'|P_{x-a}).
\end{align*}
Therefore $q$ divides $m\,e(Q'|P_{x-a})$. Since $m=q+1$, we must have that $q$ divides $e(Q'|P_{x-a})$. But from the Fundamental equality (Theorem 3.1.11 of \cite{Stichbook09}) we know that $e(Q'|P_{x-a})\leq [K(x,y):K(x)]=m$, and
this implies that $q=e(Q'|P_{x-a})$. The Fundamental equality also implies that $Q'$ and $Q$ are the only places of $K(x,y)$ lying over $P_{x-a}$. Thus the inertia degrees $f(Q'|P_{x-a})=f(Q|P_{x-a})=1$ so that $e(Q|P_{x-a})=1$. Now, since $Q$ lies over $P_{x-a}$, $\nu_Q(x-a)=1$ and $\nu_{P_{x-a}}(g(x))=0$, we have that
\begin{align*}
	m &=\nu_Q((x-a)^m/g(x))=\nu_Q((y-a)^{m}+b(y-a))\\
	&=\nu_{Q}(y-a)=e(Q|P_{y-a})\nu_{P_{y-a}}(y-a)=e(Q|P_{y-a}),
\end{align*}
 
We see that the ramification conditions required in Theorem \ref{wildramif} are satisfied with $f(X)=X-a$. Therefore 
\[\lim_{i \rightarrow \infty}\frac{g(T_i)}{[T_i:T_0]}=\infty,\]
so that $\mathcal{T}$ is a tower over $K$ of infinite genus.

\end{document}